\documentclass{amsart}
\usepackage{amsmath, amssymb}
\usepackage{amsfonts}
\usepackage[arrow,matrix,curve,cmtip,ps]{xy}

\usepackage{amsthm}

\allowdisplaybreaks

\newtheorem{theorem}{Theorem}[section]
\newtheorem{lemma}[theorem]{Lemma}
\newtheorem{proposition}[theorem]{Proposition}
\newtheorem{corollary}[theorem]{Corollary}

\newtheorem*{theorem*}{Theorem}
\theoremstyle{remark}
\newtheorem{remark}[theorem]{Remark}
\newtheorem{definition}[theorem]{Definition}

\newcommand{\thmref}[1]{Theorem~\ref{#1}}

\newcommand{\lemref}[1]{Lemma~\ref{#1}}
\newcommand{\propref}[1]{Prop\-o\-si\-tion~\ref{#1}}
\newcommand{\corref}[1]{Cor\-ol\-lary~\ref{#1}}

\numberwithin{equation}{section}

\newcommand{\C}{\mathbb{C}}

\newcommand{\go}{G^0}
\newcommand{\cc}{C_{c}}
\newcommand{\ccg}[1]{\cc(G|_{#1})}
\newcommand{\cs}{\ensuremath{C^{*}}}
\newcommand{\supp}{\operatorname{supp}}
\newcommand{\csg}[1]{\cs(G|_{#1})}
\newcommand{\ex}[1]{\operatorname{Ex}(#1)}
\newcommand{\co}{C_{0}}
\newcommand{\Prim}[1]{\operatorname{Prim}(#1)}
\renewcommand{\H}{\mathcal{H}}
\newcommand{\set}[1]{\{\,#1\,\}}
\newcommand{\Ind}{\operatorname{Ind}}

\begin{document}
\title{Classifying the Type of Principal Groupoid $\cs$-algebras}

\author{Lisa Orloff Clark }

\curraddr{Dep of Mathematical Sciences\\ Susquehanna University\\
Selinsgrove\\ PA 17870\\ USA} \email{clarklisa@susqu.edu}

\date{\today}
\begin{abstract}

Suppose $G$ is a second countable, locally compact,
Hausdorff
\newline
groupoid with a fixed left Haar system.
Let $\go/G$ denote the orbit space of $G$ and $\cs(G)$ denote
the groupoid $\cs$-algebra.  Suppose that $G$ is a principal
groupoid.  We show that $\cs(G)$ is CCR if and only if
$\go/G$ is a $T_1$ topological space, and that
$\cs(G)$ is GCR if and only if $\go/G$ is a $T_0$
topological space.   We also show that $\cs(G)$ is
a Fell Algebra if and only if $G$ is a Cartan groupoid.

\end{abstract}

\maketitle

\section{Introduction}
$\cs$-algebras can be classified as being continuous-trace, Fell Algebras, bounded trace,
 CCR (liminal), and GCR (postliminal).    These are listed in order of containment.
Recall that for separable $\cs$-algebras, an algebra is GCR if and
only if it is Type I.  Further, $\cs$-algebras that
are not GCR are very poorly behaved.
In the  case of a transformation group
$\cs$-algebra $C^*(H,X)$ (where $H$ is a group that acts continuously on the space $X$)
each of these classifications correspond to a property of
the transformation group itself.  For example, Phil Green
was able to prove in \cite{green} that a freely acting
transformation group $\cs$-algebra has continuous-trace
if and only if the action of the transformation group is proper.
In \cite{ctstrace} the authors have
generalized
Green's result to principal groupoids.  In this paper we
generalize three more such results.

In \cite{gootman}, Elliot Gootman showed  the following:
\begin{theorem}
\label{tgthm4}
Suppose $H$ and $X$ are both second countable.
Then $\cs(H,X)$ is GCR  if and only if every stability
group is GCR and the orbit space is $T_0$.
\end{theorem}

Dana Williams considered the case for CCR transformation
group $\cs$-algebras in \cite{danasthesis},
and proved the theorem below.

\begin{theorem}
\label{tgthm5}
Suppose that $H$ and $X$ are both second countable.
Suppose also that at every point of discontinuity $y$
of the map $x \mapsto S_x$, the stability group $S_y$ is amenable,
then $\cs(H,X)$ is CCR if and only if the stability
groups are CCR and the orbit space is $T_1$.
\end{theorem}
\begin{remark}
Gootman has shown that the hypothesis on $x \mapsto S_x$
in \thmref{tgthm5} is unnecessary; however, the details have not appeared.
\end{remark}
We also note that Thierry Fack proved versions of \thmref{tgthm4}
and \thmref{tgthm5} for foliation $\cs$-algebras in \cite{fack}.

 Finally, in \cite{astridsthesis}, Astrid an Huef proved :
\begin{theorem}
\label{tgthm2}
$\cs(H,X)$ is a Fell algebra if and only if (H,X) is a Cartan
$G$-space.
\end{theorem}

We generalize each of the above three
theorems to principal groupoids.  The key comes in showing that there is
a continuous injection between the orbit space of the groupoid and
the spectrum of the associated groupoid $\cs$-algebra.
In fact, when the orbit space is $T_0$, we show that these spaces are homeomorphic.

We
have also been able to further generalize the CCR and GCR results
to non-principal groupoids; however, these results will appear later.
\section{Preliminaries}

A groupoid $G$ is a small category in which every
morphism is invertible.  A principal groupoid is
a groupoid in which there is at most one morphism
between each pair of objects.    We define
maps $r$ and $s$ from $G$ to $G$ by
$r(x) = xx^{-1}$ and $s(x) = x^{-1}x$.  These are the maps Renault
calls $r$ and $d$ in \cite{renault}.
The common image of $r$ and $s$
is called the unit space which we denote $\go$.

We will only consider
second countable, locally compact, Hausdorff groupoids~$G$.
Our main results also requires $G$ to be principal;
however, we will state this condition when it is needed.  We will also assume that $G$ has a fixed
left Haar system, $\{\lambda^u\}_{u\in\go}$.

Now consider the vector space $\cc(G)$, the space of continuous functions with compact support from $G$ to the complex numbers, $\C.$  We can view
this space as a $*$-algebra by defining convolution and involution with the formulae:
\begin{align}
f*g(x) &= \int f(y)g(y^{-1}x) \ d \lambda^{r(x)}(y)\notag \\
&=\int f(xy)g(y^{-1}) \ d \lambda^{s(x)}(y)\notag
\end{align}
and
\begin{equation}
f^*(x) = \overline{f(x^{-1})}.\notag
\end{equation}

A representation of $\cc(G)$ is a $*$-homomorphism $\pi$ from
$\cc(G)$  into $B(\H)$ for some Hilbert space $\H$ that
is continuous with respect to the inductive limit topology on
$\cc(G)$ and the weak operator
topology on $B(\H)$, and that is non-degenerate in the sense that
the linear span of $\{\pi(f) \eta | f \in \cc(G), \eta \in \H\}$
is dense in $\H$.  We define the groupoid $\cs$-algebra with
the following theorem.

\begin{theorem}
\label{algebra}
For $f \in \cc(G)$, the quantity
\begin{equation}
\|f\| := \sup \{ \|\pi(f)\| \mid \pi \text{ is a
representation of } \cc(G) \}
\label{normdef}
\end{equation}
is finite and defines a $\cs$-norm on $\cc(G)$.  The completion of
$\cc(G)$ with respect
to this norm is a $\cs$-algebra, denoted $\cs(G)$.
\end{theorem}

The only real issue in proving \thmref{algebra}
comes in showing that $\|f\| < \infty $ for all
$f \in \cc(G)$.  This is a consequence of
Renault's Disintegration Theorem  \cite[Theorem~4.2]{renault2},
\cite[Theorem~3.23]{muhlysnotes}.  The motivating example of
a groupoid $\cs$-algebra is a transformation group $\cs$-algebra, $\cs(H,X)$,
defined in \cite{danasthesis} and \cite{astridsthesis}.

We define the map $\pi:G \rightarrow \go \times \go$ by
$\pi(x) = (r(x),s(x))$.  Using  $\pi$, we define an equivalence
relation on $\go$ and endow the set
of equivalence classes
with the quotient topology.  We call this topological space
the orbit space of $G$, denoted
$\go/G$.

\section{A Map from $\go/G$ to $\cs(G)^\wedge$}

Following \cite{ctstrace} and \cite[Pages~81--82]{renault},
recall that for each $u \in \go$ there
is a representation $L^u$ induced from the point mass
measure $\epsilon_u$.  When $G$ is a principal groupoid, $L^u$ acts
on $L^2(G,\lambda_u)$ so that for $f \in \cc(G)$ and
$\xi \in L^2(G,\lambda_u)$,
\begin{equation}
L^u(f)\xi(\gamma) = \int f(\gamma
\alpha)\xi(\alpha^{-1})d\lambda^u(\alpha).
\notag
\end{equation}
The following Lemma is \cite[Lemma~2.4]{ctstrace}.

\begin{lemma}
\label{welldefirred}
Suppose $G$ is a principal groupoid.
Then the representation $L^u$ is irreducible for each $u \in \go$.  Further more, if $[u]=[v]$ then $L^u$ is unitarily
equivalent to $L^v$.
\end{lemma}

 We can use this construction to define a map $\psi:\go/G
\rightarrow \cs(G)^\wedge $ where $\psi([u]) =L^u$.  As
 usual, we view $L^u$ as its unitary equivalence class in
$\cs(G)^\wedge$.
Our notation is somewhat careless.  We should
denote the image of $u$ under $\psi$ by $[L^u]$ but the preceding Lemma makes this carelessness less troubling.

Our goal is to show that for principal groupoids with $T_0$ orbit spaces, $\psi$ is a homeomorphism.  We will first show this for groupoids with $T_1$ orbit spaces and generalize this to $T_0$ orbit spaces later.  Before we deal with $\psi$, we must first determine what the representations of $\cs(G)$ look like.

Fixing $u \in \go$,
recall from \cite[Lemma~2.13]{renault} that there is a representation $M_u$ of
$\co(\go)$ on $L^2(G, \lambda_u)$  defined by
\begin{equation}
L^u(V(\phi)f)=M_u(\phi)L^u(f).
\label{mdef}
\end{equation}

\begin{proposition}
\label{ker}
Suppose that $L$ is an irreducible representation of
  $\cs(G)$, and that $M$ is the representation of $C_{0}(\go)$
  defined by $M(\phi)L(f)=L\bigl(V(\phi)f\bigr)$.  If $\ker M
=J_{F}:=\{\phi \in \co(\go) \mid \phi(x) =0
\text{ for all } x  \in
F\}$,
  then there is a $u\in\go$ such that $F=\overline{[u]}$.
\end{proposition}

Before we can prove this proposition, we need the following two
lemmas.

\begin{lemma}
  \label{lem-one}
  Let $U$ be an open subset of $\go$.  Then the ideal of $\cs(G)$
  generated by $\ccg U$ is $\overline{\ccg{[U]}}:=\ex {[U]}$.
\end{lemma}

\begin{proof}
  It suffices to see that
  \begin{align}
    E_{0}&:=\cc(G)*\ccg
    U*\cc(G)\notag\\
    &=\operatorname{span}\set{f*g*h:\text{$f,h\in\cc(G)$ and
    $g\in \ccg U$}}\label{eq:4}
  \end{align}
is dense in $\ccg{[U]}$ in the inductive limit topology.  In view
of
the Stone-Weierstrass Theorem
\cite[Theorem~7.33]{babyrudin},
since $E_0$ is self-adjoint it
suffices to show $E_0$ separates points of $G|_{[U]}$
and vanishes at no point of $G|_{[U]}$.

Because $G|_{[U]}$ is Hausdorff, this is the same as showing that
for each $\gamma \in G|_{[U]}$ and each neighborhood
$V$ of $\gamma$, there is a $F\in E_{0}$ with $\supp F\subset V$
and
$F(\gamma)\not=0$.   But if $\gamma\in G|_{[U]}$, then
$\gamma=\alpha\beta\delta$ with $\beta\in G|_{U}$, $s(\alpha) =
s(\gamma)$, and $r(\delta) = r(\gamma) $.

Now notice that
\begin{align}
  f*g*h(\gamma)&= \int_{G}
  f*g(\gamma\eta)h(\eta^{-1})\,d\lambda^{s(\gamma)} (\eta) \notag\\
&=\int_{G}\int_{G} f(\omega)g(\omega^{-1}\gamma\eta) h(\eta^{-1})
  \,d\lambda^{r(\gamma)} (\omega)
\,d\lambda^{s(\gamma)}(\eta).\notag\\
&=\int_{G}\int_{G} f(\omega)g(\omega^{-1}\gamma\eta^{-1}) h(\eta)
  \,d\lambda^{r(\gamma)} (\omega)
\,d\lambda_{s(\gamma)}(\eta).\notag
\end{align}

We can choose neighborhoods $V_{1}$, $V_{2}$ and $V_{3}$ of
$\alpha$,
$\beta$ and $\delta$, respectively, such that $V_{1}V_{2}V_{3}
\subset V$.
Notice from the integral above that if $\gamma \in \supp(f*g*h)$
then there exists $\omega \in \supp f$,
$\eta \in \supp h$ so that $\omega^{-1}\gamma \eta^{-1}
\in \supp g$.  Since $\gamma = \omega (\omega^{-1}\gamma)
\eta^{-1}
\eta$,
we see that $\supp(f*g*h)\subset (\supp f)(\supp g)(\supp h)$, so we
have $\supp
(f*g*h)\subset V$ provided $\supp f\subset V_{1}$, $\supp g\subset
V_{2}$ and
$\supp h\subset V_{3}$.  Thus it suffices to take
non-negative functions $f,h\in\cc(G)$ and $g\in \ccg U$ with the
appropriate supports and $f(\alpha)=g(\beta)=h(\delta)=1$ and
$F=f*g*h$.
\end{proof}

\begin{lemma}
  \label{lem-two}
  Suppose that $L$ is a non-degenerate representation of
  $\cs(G,\lambda)$, and that $M$ is the representation of
$C_{0}(\go)$
  defined by $M(\phi)L(f)=L\bigl(V(\phi)f\bigr)$.   Then $\ker
  M=J_{F}$ for a closed, $G$-invariant set $F\subset \go$.
\end{lemma}

\begin{proof}
  We know $\ker M=J_{F}$ for closed subset $F$ of $\go$.
  Let $U:= \go\setminus F$.  It will suffice to see that $U$ is
  $G$-invariant; that is, $U=[U]$.

If $f\in\ccg U$, then $K=\supp f$ is a
  compact subset of $G|_{U}$.  Thus $C=r(K)$ is a compact subset
of
  $U$.  Therefore we can choose $\phi\in\cc(U)$ such that $\phi(u)
=1$ for
  all $u\in C$.  Then $V(\phi)f=f$.  Since $\phi$ vanishes on
$F$, $M(\phi)L(f) =L(V(\phi)f) = 0$.  So
  $f\in\ker L$, and we
  have shown that
  \begin{equation}
    \label{eq:1}
    \ccg U\subset \ker L.
  \end{equation}
Lemma~\ref{lem-one} implies that $\ccg{[U]}\subset \ker L$. If
$[U]\not=U$, then there is a $\phi\in \cc(\go)$ such that $\supp
\phi\subset [U]$ and $\phi$ is not identically zero on $F$.  Since
$V(\phi)f \in \ccg {[U]}$ for all $f\in \cc(G)$, it follows that
$V(\phi)f \in \ker L$. Therefore $M(\phi)=0$, which contradicts
$\ker M=J_{F}$.
\end{proof}

\begin{proof}[Proof of \propref{ker}]
Since $\go/G$ is a second countable Baire space,
we know from \cite[Lemma on page 222
preceding Corollary~19]{green} every irreducible
closed set must be a point closure.
Lemma~\ref{lem-two} tells us that $\ker M=J_{F}$ where $F$ is a
closed
$G$-invariant subset of $\go$.  Thus the image of $F$ in $\go/G$ is
closed.
Suppose F is not an orbit closure. Then F is not irreducible.  That
is F can be written as the
union $C_{1}\cup C_{2}$ where each $C_{i}$ is a closed
$G$-invariant
set such that $F\not\subset C_{i}$.  In particular,
$C_{i}\cap F\not=\emptyset$ for $i=1$ or $i=2$.

Let $U_{i}$ be the $G$-invariant open set $ \go\setminus C_{i}$.
 Since $\ex{U_{1}}\cap \ex{U_{2}} = \ex{U_{1}}\ex{U_{2}}$, it
follows
from \cite[Lemma~2.10]{ctstraceIII} that
\begin{equation*}
  \ccg{U_{1}}\ccg{U_{2}}
\end{equation*}
is dense in $\cs(G|_{U_{1}})\cap \cs(G|_{U_{2}})$.  On the other
hand
\begin{align*}
  \ccg{U_{1}}\ccg{U_{2}} &\subset \ccg{U_{1}\cap U_{2}} =
  \ccg{\go\setminus (C_{1}\cup C_{2})}\\
  & = \ccg{\go\setminus
  F}=\ccg{U}.
\end{align*}
Thus, \eqref{eq:1} implies that
\begin{equation*}
  \ex{U_{1}}\cap \ex{U_{2}} \subset \ker L.
\end{equation*}
Since $L$ is irreducible, $\ker L$ is prime.  Thus
\begin{equation*}
  \ex{U_{i}}\subset \ker L\quad\text{for some $i=1,2$}.
\end{equation*}
We may as well assume that $i=1$.  Since $U_{1}\cap F\not
=\emptyset$
(otherwise, we'd have $F$ in $C_{1}$),
we can choose $\phi\in\cc^{+}(\go)$ such that $\supp\phi\subset
U_{1}$
and $\phi|_{F}\not=0$.  If $f\in\cc(G)$, we know
    $$V(\phi)f(\gamma) = \phi(r\gamma)f(\gamma)$$
thus $r(\gamma) \in U_1$ and because $U_1$ is invariant, $s(\gamma)
\in U_1$ also.
This means that $V(\phi)f$ is in $\cc(G|_{U_{1}})$.  Thus
$V(\phi)f\in\ker
L$ for all $f\in \cc(G)$.  It follows that $M(\phi)=0$.  But this
contradicts $\phi|_{F}\not=0$.  Thus $F$ must be an orbit
closure as claimed.
\end{proof}

\begin{corollary}
\label{factors}
Every irreducible representation  of
$\cs(G)$
factors through $\csg{\overline{[u]}}$ for some $u \in \go$.
\end{corollary}

\begin{proof}
Suppose $L$ is an irreducible representation and M is the
associated
representation satisfying (\ref{mdef}).
We know $\ker M=J_{F}$ and that $F = \overline{[u]}$ by
\propref{ker}.
 Let $U:= \go \setminus F$.
We must  show that $\ex U \subset \ker L$ by \cite[Lemma~2.10]{ctstraceIII}.  It
suffices to show $\cc(G|_U) \subset \ker L$.
We will do this as we did in the proof of \lemref{lem-two}.
If $f\in\ccg U$, then $K=\supp f$ is a
  compact subset of $G|_{U}$.  Thus $C=r(K)$ is a compact subset
of
  $U$.  Therefore we can choose $\phi\in\cc(U)$ such that $\phi(u)
=1$ for
  all $u\in C$.  Then $V(\phi)f=f$.  Since $\phi$ vanishes on
$F$, $M(\phi)L(f) =L(V(\phi)f) = 0$.  So
  $f\in\ker L$, and we
  have shown that
    $\ccg U\subset \ker L.$
\end{proof}

We now have all the pieces needed to show that for principal
groupoids,  the map $\psi$ is a continuous open injection.  Further, if the orbit space is $T_1$, then $\psi$ is a homeomorphism.

\begin{proposition}
\label{cio}
Suppose $G$ is a principal groupoid.
Then the map $\psi$ defined above is a continuous, open, injection.
\end{proposition}

\begin{proof}
We know that $\psi$ is a continuous injection by \cite[Propostion~2.5]{ctstrace}.

We will show $\psi$ is an open map using the criteria from
\cite[Proposition~II.13.2]{imprimitivity}.
Let $L^{u_n} \rightarrow L^u$ be a convergent net in $\cs(G)^\wedge
$.  Thus $M_{u_n} \rightarrow M_u$ in $\co(G^0)^\wedge$.
Each $M_{u_n}$ corresponds to a closed subset, namely
$\overline{[u_n]}$.  By \cite[Lemma~2.4]{danasthesis},
we may pass to a subnet and relabel if necessary and find
$v_n \in [u_n]$ so  $v_n \rightarrow u$.
Therefore $\psi$ is open.
\end{proof}

\begin{remark}  We will eventually weaken the hypothesis of \propref{homeo} and
require only that $G$ be a principal groupoid and $\go/G$ be $T_0$.
\end{remark}

\begin{proposition}
\label{homeo}
Suppose $G$ is a principal groupoid in which orbits are closed.
Then the map $\psi$ defined above in a homeomorphism.
\end{proposition}

\begin{proof}
All that is left to show is that $\psi$ is surjective.
Let $L$ be any irreducible representation of $\cs(G)$.
Since orbits are closed,  we know that $L$ is lifted from a
representation on $\csg{[u]}$ from
\corref{factors}.
The representation $L^u$ is also a representation on $\csg{[u]}$.
Since
$\csg{[u]}$ is a transitive groupoid,
and $G$ is principal,  \cite[Lemma~2.4]{equivalence} tells us that
$\csg{[u]} \cong K(H)$.  However, the compact
operators have only one irreducible representation.  Therefore $L^u
\cong L$.
\end{proof}


\section{CCR Groupoid $\cs$-algebras}

In order to prove the theorem below, a generalization of Williams' \thmref{tgthm5}, we only use the property of \propref{cio} that $\psi$ is a continuous injection.

\begin{theorem}
\label{ccrthm}
Let $G$ be a principal groupoid.  Then $G$ is CCR if and
only if $\go$ is~$T_1$.
\end{theorem}

\begin{proof}
Suppose $C^*(G)$ is CCR.  This implies that points of the spectrum,
$C^*(G)^\wedge$, are closed.  We know the map,
$$\psi:G^0/G \rightarrow C^*(G)^\wedge$$
where $\psi([u]) = L^u$ is a continuous injection by
\propref{cio}.  Thus the inverse image of a point
of
the spectrum is one orbit which must also be closed.

Now suppose that the orbit space is $T_1$.
Suppose $L$ is a representation of $\cs(G)$.
We know from\corref{factors} that $L$ factors through
$\cs(G|_{\overline{[u]}})=\cs(G|_{[u]})$ for some
$u \in \go$. But $\cs(G|_{[u]})$ is a transitive groupoid thus
\begin{equation}
\cs(G|_{[u]})\cong \cs(G^u_u) \otimes K\notag
\end{equation}
by \cite[Theorem~3.1]{equivalence}. This is CCR because
we are assuming $G$ is a principal groupoid. This means that
$L$  is lifted from a representation  of a CCR $\cs$-algebra
making $L$ a representation onto the compact operators.
That is, $\cs(G)$ is CCR.
\end{proof}

\begin{corollary}
\label{ccrhomeo}
If $G$ is a principal groupoid and $\cs(G)$ is CCR then $\psi$ is a
homeomorphism.
\end{corollary}

\begin{proof}
This is immediate from \thmref{ccrthm} and
\propref{homeo}.
\end{proof}

\section{GCR $\cs$-algebras}

We can weaken the conditions in Proposition~\ref{homeo}
and show that, for principal groupoids,
$\psi$ is a homeomorphism when $\go/G$ is a $T_0$ space.
In doing this, we actually describe the ideal
structure of the associated groupoid
$\cs$-algebra.
We will also prove a generalization of Gootman's \thmref{tgthm4}
for principal groupoids  that says
$\cs(G)$ is GCR if and only if $\go/G$ is $T_0$.

We know that for principal groupoids
$\psi$ is a continuous, injective, open map from \propref{cio}.
Therefore to show $\psi$ is a homeomorphism,
we must show that $\psi$ is onto.  What we will
do is show that when we require the orbit
space to be $T_0$ rather than $T_1$, we can show that every irreducible
representation of
$\cs(G)$ is lifted from a representation of $\cs(G|_C)$ where C is a
Hausdorff subset of $\go/G$.  This will suffice.

We will begin Proposition~\ref{t0homeo} below by assuming that
$\go/G$ is $T_0$.  We will also show that the orbit equivalence relation
R on $\go$ is an $F_\sigma$
subset of $\go \times \go$.  When this is the case, Arlan Ramsay
proved in \cite[Theorem~2.1]{ramsay}
that there is a list of 14 different properties that are each
equivalent to saying that $\go/G$ is $T_0$.
Some of these equivalent
properties include:  (1) each orbit is locally closed, (2) $\go/G$ is
almost Hausdorff, and (3) $\go/G$ is a standard Borel space.  We will use
property (2) in our proof.  The idea for this proof comes from \cite[Lemma~2.3]{orbits}.

\begin{proposition}
\label{t0homeo}
Suppose $G$ is groupoid.  If $\go/G$ is $T_0$ then there is an
ordinal $\gamma$ and
ideals $\{I_\alpha :\alpha \leq \gamma\}$ such that
\begin{enumerate}
\item[(i)]  $\alpha < \beta$ implies that $I_\alpha \subset I_\beta$,

\item[(ii)]  $I_0 = 0$ and $I_\gamma = \cs(G)$,

\item[(iii)]  if $\delta$ is a limit ordinal, then
$I_\delta$ is the ideal generated by $\{I_\alpha\}_{\alpha<\delta}$,

\item[(iv)]  if $\alpha$ is not a limit ordinal, then
$I_{\alpha} / I_{\alpha-1} \cong \cs(G|_{U_{\alpha}
\backslash U_{\alpha-1}})$ where $U_\alpha$ is a saturated subset of $G$ and
each space $U_{\alpha +1} \backslash U_\alpha$ is Hausdorff and

\item[(v)]   if $L$ is an irreducible representation of $\cs(G)$, then $L$
is the canonical extension of an irreducible
 representation of $\cs(G|_{U_{\alpha} \backslash U_{\alpha-1}})$.
\end{enumerate}

Also, if $G$ is a principal groupoid, then the map $\psi$ defined
above is a homeomorphism from $\go/G$ into
 $\cs(G)^\wedge$.

\end{proposition}

\begin{remark}
The $\cs$-algebra $\cs(G|_{U_{\alpha} \backslash U_
{\alpha-1}})$ is actually the quotient of
$\cs(G|_{U_\alpha})$ by $\cs(G|_{U_{\alpha-1}}).$
\end{remark}

\begin{proof}
First we will show that the orbit equivalence relation R on $\go$ is an
$F_\sigma$ subset of $\go \times \go$.  To show that R is an $F_\sigma$ set,
we must show it is a countable union of closed sets of $\go \times \go$.
Notice that $G$ is $\sigma$-compact and that $R = \pi(G)$ where
$\pi(\gamma) = (r(\gamma), s(\gamma))$.  Therefore $R$ is an $F_\sigma$
subset because $\pi$ is continuous.

Now from \cite[Theorem~2.1]{ramsay}, we know that $\go/G$ is almost Hausdorff.
Therefore, the discussion on page 125 of \cite{Glimm} gives us an ordinal
$\gamma$ and open subsets $\{U_\alpha : \alpha \leq\gamma\}$ of
$\go/G$ such that

(a)  $\alpha < \beta  $ implies that $U_\alpha \subset U_\beta$

(b)  $\alpha < \gamma$ implies that $U_{\alpha} \backslash
U_{\alpha-1}$ is a dense Hausdorff subspace in the relative
topology.

(c)  if $\delta$ is a limit ordinal, then
$$U_\delta = \underset{\alpha<\delta}{\bigcup} U_\alpha $$

(d)  $U_0 = \emptyset$ and $U_\gamma = \go/G$

In the sequel, we will abuse notation and consider each
$U_\alpha$ as an open invariant
subset of $\go$.  Thus from \propref{idealbijection}
each $U_\alpha$corresponds to an ideal,
$\cs(G|_{U_\alpha})$ of $\cs(G)$ which we will call $I_\alpha$.
Now properties (i), (ii), and (iii) follow immediately.
Property (iv) follows immediately from the short exact sequence
\begin{equation}0 \longrightarrow\cs(U|_{\alpha-1})
\longrightarrow\cs(G|{U_\alpha})
\longrightarrow \cs(G|_{U_\alpha\setminus U_{\alpha-1}})
\longrightarrow
0\notag
\end{equation} of
\cite[Lemma~2.10]{ctstraceIII}.

Now we must show (v).  Suppose $L$ is an irreducible representation of
$\cs(G)$. Since $L$ is an non-degenerate irreducible representation, the
restriction of $L$ to an ideal gives us an irreducible representation of
the ideal.  Define the set
\begin{equation}
S=\{\lambda \mid L(I_\lambda)
\neq 0\}.\notag
\end{equation}
Since $S$ is a set of ordinals, it has a
smallest element.  Let $\alpha$ be the smallest element
of $S$.  We know that $\alpha$ is not a limit
ordinal because property of (iii).  Therefore $\alpha-1$ exists and we
have
\begin{equation}
L(I_\alpha)\neq 0 \text{ and } L(I_{\alpha-1}) =
0.\notag
\end{equation}
Therefore, $L$ is the canonical extension of a
representation of $I_\alpha/I_{\alpha-1}$ as needed.

Suppose $G$ is a
principal groupoid. We know that $\psi$ is continuous, open,  and
injective from
\propref{cio}.
Thus, to show $\psi$ is a homeomorphism,  we need only show that $\psi$ is onto.
In this proof, we need to be careful and define the following
representations.  Let $\Ind (G, u)$ be the representation
$L^u$ on $\cs(G)$ and let $\Ind(G_{U_\alpha}, u)$ be the representation $L^u$ as a
representation of $\cs(G|_{U_\alpha})$ for some $u \in U_\alpha$.

Now let $L$ be any representation of $\cs(G)$.  Our goal is to
show that $L$ is equivalent to $L^u=\Ind (G, u)$ for some $u \in \go$.  We
know from part (v) that $L$ is the canonical extension of a representation
$L'$ of $I_\alpha/I_{\alpha-1}= \cs(G|_{U_\alpha\setminus
U_{\alpha-1}})$.  We also know that $U_\alpha\setminus U_{\alpha-1}$ is
Hausdorff which means that $L'$ is equivalent to
$\Ind(G|_{U_\alpha},u)$ for some $u \in U_\alpha$.  It suffices to show
that the canonical extension of $\Ind(G|_{U_\alpha},u)$ to
$\cs(G)$ must be equal to
$\Ind(G,u)$.  Notice that the spaces each of these representation act upon
are the same.  The representation $\Ind(G|_{U_\alpha},u)$  extends to  a
representation  $\overline{\Ind}(G|_{U_\alpha},u)$  on all
of $\cs(G)$.
Notice that for $f \in \cc(G)$, $g \in L^2(G_u,\lambda_u)$,
$x \in G_u$
we have
\begin{align}
\overline{\Ind}(G|_{U_\alpha},u)(f)(\Ind(G|_{U_\alpha},u)(g))
\xi &=
\overline{\Ind}(G|_{U_\alpha},u)(f*g)\xi\notag\\
&=\Ind(G,u)(f*g)\xi.
\notag
\end{align}
Thus, $\Ind(G,u)$ is the canonical extension of $\Ind(G|_{U_\alpha},u)$
as needed.
\end{proof}

We now have more than enough to prove the following theorem.

\begin{theorem}
\label{gcrprinc}
Suppose $G$ is a principal groupoid.  Then $\cs(G)$ is GCR if and
only if $\go/G$ is $T_0$.
\end{theorem}

\begin{proof}
Suppose $\cs(G)$ is GCR.  Then the spectrum of $\cs(G)$ is $T_0$.
From \lemref{homeo}, we know there is a continuous injection from
the orbit space into the spectrum. Therefore, the orbit space must
also be $T_0$.

Now suppose we know $\go/G$ is $T_0$.  From \propref{t0homeo}, we
know that every irreducible representation L of
$\cs(G)$ is the canonical extension of a representation of
$\cs(G_{U_\alpha \setminus U_{\alpha-1}})$ where
$U_\alpha \setminus U_{\alpha-1}$ is Hausdorff.
Thus $\cs(G_{U_\alpha \setminus U_{\alpha-1}})$ is CCR
by \thmref{ccrthm}.
Therefore, the image of $L$ contains the compact operators
and $\cs(G)$ is GCR.
\end{proof}

\section{ideals}
We know that for an open saturated subset $U$ of $G^0$,
$C^*(G|_U)$ is an ideal in $C^*(G)$.  When $G$ is
principal and $\cs(G)$ is GCR, all the ideals of
$C^*(G)$  are of this form.

\begin{proposition}
\label{idealbijection}
Suppose $G$ is a principal groupoid and $\cs(G)$ is GCR.  Then the
map $U \mapsto \ex U \cong C^*(G|_U)$ from the
collection of open saturated subsets of $G^0$ to the ideals of
$C^*(G)$ is a bijection.
\end{proposition}

\begin{proof}
Recall that if $\cs(G)$ is GCR, $\cs(G)^\wedge \cong
\Prim{\cs(G)}$. We also know that there is a natural
correspondence between open subsets of $\Prim {\cs(G)}$ and ideals
of $\cs(G)$.  Thus in order to show that Ex is a
bijection, it suffices to show
\begin{equation}
C^*(G|_U) \cong \underset{v \notin U}{\bigcap} \ker L^v.
\notag
\end{equation}

Notice that
\begin{equation}
C^*(G|_U) = \bigcap\{\ker L^v : L^v(C^*(G|_U))=0\}.
\notag
\end{equation}

It follows from the definition of $L^v$ that
$\text{if } v\in U \text{, } L^v(\cc(G|_U)) \not= 0 \text{  and
} \text{if } v \notin U \text{, } L^v(C^*(G|_U)) = 0$.  Therefore
\begin{equation}
    C^*(G|_U) = \underset{v \notin U}{\bigcap} \ker L^v\notag
\end{equation}
as needed.
\end{proof}

\section{Fell Algebras}

Finally, we generalize an Huef's
\thmref{tgthm2}.
Many of the results involving Cartan $G$-spaces that an Huef used to
prove
\ref{tgthm2} came from \cite{palais}.   Thus
we first must generalize some of Palais' work for Cartan G-spaces.
This process leads us to some interesting
results in their own right.

\begin{definition}
\label{wandering}
A subset, $N$ of $G^0$ is wandering if and only if the set
\begin{equation}
G|_N = \pi^{-1}(N,N) = \{\gamma \in G \mid s(g)\in N \text{ and } r(g) \in N\}\notag
\end{equation}
is relatively compact.
\end{definition}

\begin{lemma}
\label{properwandering}
A groupoid $G$  is proper if and only if every compact subset
of $G^0$ is wandering.
\end{lemma}

\begin{proof}
Suppose $G$ is proper so that by definition $\pi$ is a proper map.
That is, the inverse image of a
compact set is compact.  Let $K$ be a compact subset of $G^0$.
By assumption  $\pi^{-1}(K,K)$  is compact; thus $K$ is wandering.

Now suppose that every compact subset of $G^0$ is wandering.  Let
$L$ be a compact subset of $G^0 \times G^0$.
We must show $\pi^{-1}(L)$ is compact. Note that $L \subset W
\times W$ where $W$ is a compact subset of $G^0$.

Thus,
\begin{equation}
\pi^{-1}(L) \subset \pi^{-1}(W,W)\notag
\end{equation}
which is compact. Thus
$\pi^{-1}(L)$ is a closed subset of a compact set.
Therefore $\pi^{-1}(L)$ is compact.
\end{proof}

\begin{definition}
\label{cartan}
We call a groupoid $G$ a Cartan groupoid if and only if for every
$x
\in G^0$, $x$ has a wandering neighborhood.
\end{definition}

It is not difficult to show that a transformation group is a Cartan
$G$-space
if and only if the
associated
transformation group groupoid is a
Cartan groupoid.

\begin{lemma}
\label{orbitsclosed}
If G is a Cartan groupoid, then for each $u\in \go$, $[u]$ is closed
in $G^0$.
\end{lemma}

\begin{proof}
Let $u \in G^0$.  Let $v$ be a limit point of $[u]$ in $G^0$.
Because $G$ is a Cartan groupoid, $v$ has a  wandering
 neighborhood, $U$. We will assume that $U$ is closed.
Thus, we can find a sequence of elements,
 $\{v_n \}$ in $U$ that converge to $v$
where each $v_n \in [u]$.  There also exists a sequence of
elements, $\{\gamma_n \} \subset G$ such that for each $n$,
 $s(\gamma_n ) = v_n $ and
$r(\gamma_n ) = u$. Now choose one of the  $\{\gamma_n \}$, call
it   $\gamma_{n_0}$.  Notice that
$r(\gamma_{n_0} ^{-1}) = v_{n_0}$ and $s(\gamma_{n_0}^{-1}
)=u$. Thus $ \gamma_{n_0}^{-1} \gamma_n \in G|_U$ which is
compact because
it is relatively compact and closed .  Thus we can pass to
 a subsequence, relabel, and assume $\{\gamma_n\}$ converges to
$\gamma$.  Since $r$ and $s$ are continuous, $r(\gamma)=u$ and
$s(\gamma)=v$.
Thus $v \in [u]$.
\end{proof}

Clearly, if $G$ is proper, by \lemref{properwandering} we see that
$G$ is a Cartan groupoid.
 We will prove a partial converse of this but first we need the
following lemma.

\begin{lemma}
\label{propernets}
A groupoid $G$ is proper if and only if every sequence,
$\{\gamma_n\}
\in G$ such that $\{\pi(\gamma_n)\}$ converges  has a convergent
subsequence.
\end{lemma}

\begin{proof}
Suppose that $G$ is proper.  Let $\{\gamma_n\}$ be a sequence where
$\{\pi(\gamma_n)\}$ converges to $(u,v)$.  Now, let $K$ be
 a compact neighborhood of $(u,v)$.  Thus
$\{\pi(\gamma_n)\}$ is eventually
inside of $K$.  Since $\pi^{-1}(K)$ is compact,
there is a subsequence,
$\{\gamma_{n_k}\}$  that converges to $\gamma$ as needed.

Now suppose for every $\{\gamma_n\} \in G$ such that
$\pi(\gamma_n)$ converges to $(u,v)$, $\{\gamma_n\}$ has a
convergent
 subsequence  $\{\gamma_{n_k}\}$ where  $\{\gamma_{n_k}\}$
converges to
$\gamma$.  Let $K$ be a compact subset of $\go \times \go$.
We must show $\pi^{-1}(K)$ is
compact.  Let $\{\gamma_n\} \subset \pi^{-1}(K)$.  It suffices to
show $\{\gamma_n\}$ has a convergent subsequence.
Since $\{\pi(\gamma_n)\} \subset K$,
$\{\pi(\gamma_n)\}$ has a convergent subsequence in $K$, call it
$\{\pi(\gamma_{n_k})\}$ where
$\{\pi(\gamma_{n_k})\} \rightarrow (u,v)$.  So, by assumption, we
can find
a subsequence and relabel so that $\{\gamma_{n_{k}}\}$ converges to
$\gamma \in \pi^{-1}(K)$.
\end{proof}

\begin{lemma}
\label{hausproper}

A groupoid $G$ is proper if and only if $G$ is Cartan and $G^0/G$
is Hausdorff.
\end{lemma}

\begin{proof}
Suppose $G$ is Cartan and $G^0/G$ is Hausdorff.  Let
$\{\gamma_n\}$ be a sequence in $G$  such
that $\{\pi(\gamma_n)\}$ converges to $(u,v)$.  By
\lemref{propernets}, we must show that there
 exists a convergent subsequence of $\{\gamma_n\}$ that converges
to
$\gamma$.

Because the quotient map is continuous,
\begin{equation}
[r(\gamma_n)] \rightarrow [u] \text{ and }[s(\gamma_n)]
\rightarrow [v]
\notag
\end{equation}
in $\go/G$.  Since the orbit space is Hausdorff, and for each $n$
\begin{equation}
[r(\gamma_n)] = [s(\gamma_n)],
\notag
\end{equation}
we must have $[u]=[v]$.  Thus there exist $\gamma \in G$ so that
$r(\gamma)=u$ and $s(\gamma) = v$.  Which also means that
\begin{equation}
r(\gamma_n) \rightarrow r(\gamma) \text{ and }s(\gamma_n)
\rightarrow s(\gamma).
\notag
\end{equation}
That is,
\begin{equation}
\pi(\gamma_n) \rightarrow \pi(\gamma)=(u,v).
\notag
\end{equation}
Since $r$ is open, we can pass to a subsequence, relabel, and find
$\eta_n \rightarrow \gamma$ with $r(\eta_n) = r(\gamma_n)$.
Then $\eta^{-1}_n\gamma_n$ makes sense and
$\pi(\eta^{-1}_n\gamma_n) \rightarrow (v,v)$.
By taking a wandering neighborhood $U$ of $v$, we can pass to a
subsequence,
relabel, and assume that $\eta^{-1}_n\gamma_n \rightarrow \beta$
with $\beta \in G|_{\{v\}}$.  But then
$\gamma_n \rightarrow \gamma\beta$ as needed.

Now suppose $G$ is proper.  Since $G$ is locally
compact, \lemref{properwandering} tells us that $G$
is Cartan.  We must show that $\go/G$ is Hausdorff.
It suffices to show that limits of convergent nets are unique.

Suppose  $\{x_n\} \in \go$ and
\begin{equation}
[x_n] \rightarrow [u] \text{ and } [x_n] \rightarrow [v].
\notag
\end{equation}
Notice that the quotient map
\begin{equation}
q: \go \rightarrow \go/G\notag
\end{equation}
is open. This is true because $q(U) =s(r^{-1}(U))$ for any
open set $U \in \go$ and $r$ and $s$ are continuous
and open.   Thus using \cite[Propostion~2.13.2]{imprimitivity},
we can pass to a subnet, relabel, and assume that
$x_n$ converges to $x$ in $\go$ and that there are
$\{v_n\} \subset \go$ such that $[v_n] = [x_n]$ with
$v_n$ converging to some $v$.  Similarly, we can find
$\{u_n\}\subset \go$ such that $[u_n] = [x_n] = [v_n]$.

Let $\gamma_n \in G$ be such that $r(\gamma_n) = u_n$ and
$s(\gamma_n) = v_n$.  If $K$ is a compact neighborhood of $u$ and $v$, then $\{\gamma_n\}$ is eventually in the compact set
$\pi^{-1}(K,K)$.  Thus we can pass to a subnet,
relabel, and assume that $\gamma_n$ converges to $\gamma$ in
$G$.  But then $(\gamma)=u$ and $s(\gamma)=v$. That is
$[u]=[v]$.
\end{proof}

Because of the correspondence between open saturated subsets and ideals, saturated sets give us a key to the structure of $\cs(G)$.  For Cartan groupoids, we can take  the
saturation of wandering neighborhoods and see that in addition to getting a saturated set,  some of the
useful properties of wandering neighborhoods are preserved.

\begin{lemma}
\label{homeoorbit}
Suppose $G$ is a principal Cartan groupoid and $U$ is an open
wandering neighborhood.  Let $V := [U]$ be
the saturation of $U$. Then $V/G|_V$ and $U/G|_U$ are homeomorphic.
\end{lemma}

\begin{proof}

Suppose that
$$q_U :U \rightarrow U/G|_U$$
and
$$q_V :V \rightarrow V/G|_V$$
are the corresponding quotient maps for the orbit spaces for $G|_U$
and $G|_V$.  Now consider the map
$$f:U/G|_U \rightarrow V/G|_V$$ so that
$$f(q_U(x)) = q_V(x)$$
for $x \in U$.  We will show $f$ is a homeomorphism.  Clearly, $f$
is
well defined.

Suppose $$q_V(x_1) = q_V(x_2) \text{ where } x_1, x_2 \in U$$
This means there exist $\gamma \in G|_V$ so that $r(\gamma) =
x_1$ and $s(\gamma) = x_2$.
Since we know $x_1$ and $x_2$ are in $U$, $\gamma \in G|_U$.
Therefore
$$q_U(x_1) = q_U(x_2)$$
and $f$ is injective.

Now let $q_V(y)\in V/G|_V$. Since $y \in V$ and $V=[U]$, $y$ is
in the orbit of $x$ for some $x \in U$.
This means that $q_V(y) = q_V(x)=f(q_U(x))$ and f is
surjective.

Suppose that $\{q_U(x_n)\}$ converges to $q_U(x)$.  We must show
that $\{q_V(x_n)\}$ converges to $q_V(x)$.
Suppose the contrary.  Thus we can find a neighborhood, $W$, of
$q_V(x)$ for which there is a subsequence which we relabel
and assume $\{q_V(x_n)\} \notin W$ for all $n$.  Because
$\{q_U(x_n)\}$ converges to $q_U(x)$, and $q_U$ is an open map,
it follows from \cite[Proposition 2.13.2]{imprimitivity} that we
can find a sequence $\{y_n\}$ and a subsequence of $\{x_n\}$ and
relabel so that $y_n \rightarrow x$ and $[y_n]=[x_n]$ in $U$.
Therefore $q_V(y_n) = q_V(x_n)$ for all $n$ and,
since  $q_V$ is continuous, $\{q_V(x_n)\}$ converges to $q_V(x)$.
This is a contradiction; thus $f$ is continuous.

Suppose $q_V(u_n) \rightarrow q_V(u)$ where we can suppose that
each
$u_n$ as well as each $u$ belong to $U$.
Since $q_V$ is open, we can pass to a subsequence, relabel, and
assume
that there are $v_n$ in $V$ such that
$q_V(v_n)=q_V(u_n)$ and $v_n \rightarrow u$.  Since $U$ is open,
we eventually have each $v_n \in U$.
Since $q_U$ is continuous, for large $n$, $q_U(v_n) \rightarrow
q_U(u)$.  It follows from \cite[Proposition~II.13.2]{imprimitivity}
that $f$ is open.

\end{proof}

\begin{lemma}
\label{saturatedproper}
Suppose $V$ is the saturation of an open wandering set, then $G|_V$
is proper.
\end{lemma}

\begin{proof}
Because $G$ is a Cartan groupoid, $G|_V$ is also a Cartan groupoid.
Thus, to show that $G|_V$ is proper, it suffices
to show that the orbit space, $V/G|_V$, is Hausdorff.  From
\lemref{properwandering}, we know that $G|_U$ is proper,
thus by \lemref{hausproper},  $U/G|_U$ is Hausdorff.  But
\lemref{homeoorbit} tells us that $U/G|_U  \cong V/G|_V$.
Therefore $V/G|_V$ is also Hausdorff.
\end{proof}

With this newly defined structure of a Cartan groupoid, we have the machinery to generalize \thmref{tgthm2}.

\begin{theorem}
\label{fellthm}
Suppose $G$ is a principal groupoid.  Then $G$ is a Cartan Groupoid
if and only if $A=\cs(G)$ is a Fell algebra.
\end{theorem}

\begin{proof}
Suppose $G$ is a Cartan groupoid. We must show that for every
irreducible representation, $\pi$ of $A$, $\pi$ is a Fell
point of $\hat{A}$.
Let $x \in G^0$ and $U$ be an open wandering neighborhood of $x$.
Let $V$ be the saturation of $U$ which is also open.

Since $G$ is a Cartan groupoid, the orbits of $G$ are closed by
\lemref{orbitsclosed}.  Therefore $G^0/G \cong \widehat{A}$
 by \propref{homeo}.
Let $\pi$ be the representation of $A$ that corresponds to $[x]$.

Since $V$ is a saturated open subset of $G$, \cite[Lemma~2.10]{ctstraceIII} tells us
$C^*(G|_V)$ is an ideal in $A$.
Thus $\pi$ is an irreducible representation of $C^*(G|_V)$.  Also,
from \lemref{saturatedproper}, we know
that $G|_V$ is a principal proper groupoid thus
\cite[Theorem~2.3]{ctstrace}, tells us that the ideal,
$C^*(G|_V)$ has continuous-trace.
We know continuous-trace $\cs$-algebras are
Fell algebras, thus $\pi$ is a Fell point of the open subset
$C^*(G|_V)^\wedge$ of $\hat{A}$ which means $\pi$ is Fell point of
$\hat{A}$ also.

Now suppose $A$ is a Fell algebra.  Let $x \in G^0$.  We must show
$x$ has a wandering neighborhood.

Since $A$ is CCR, $G^0/G \cong \widehat{A}$ by \corref{ccrhomeo}.

Let $\pi_x$ be the representation corresponding to $[x]$.  Since
$\pi_x$ is a Fell point, from \cite[Corollary~3.4]{archbold}
we know it has an open Hausdorff neighborhood in $\widehat{A}$.
This neighborhood is of the form $\widehat{J}$ where $J$ is
an ideal of $A$.  We also know from \lemref{idealbijection} that
\begin{equation}
J \cong C^*(G|_V)\notag
\end{equation}
for some open, saturated subset $V$ of
$G^0$.  Notice that $x \in V$.

Since $J$ has Hausdorff spectrum and is a Fell algebra,
$J$ has
continuous-trace.  Therefore by \cite[Theorem~2.3]{ctstrace},
$G|_V$ is proper. Thus, we know from \lemref{properwandering} that
every compact subset of $V$ is wandering.

Let $N$ be a compact neighborhood of $x$ in $V$.
Therefore $N$ is
wandering neighborhood of $x$ in $G^0$.
\end{proof}

The proof of the
following corollary is trivial in the transformation group case
however requires much of the
machinery established thus far to prove it in the groupoid case.

\begin{corollary}
\label{cor1}
Suppose $G$ is a principal groupoid.  If $x \in \go$ has a
wandering neighborhood and $y \in [x]$,
then $y$ has a wandering neighborhood.
\end{corollary}

\begin{proof}
Let $U$ be an open wandering neighborhood of $x$.  We know that
$G|_{[U]}$ is proper. Therefore
$\csg{[U]}$ has continuous-trace  which means it is a Fell algebra.
Thus by \thmref{fellthm}, $G|_U$ is a Cartan groupoid.  So we know
every element of $[U]$
has a wandering neighborhood in $[U]$; therefore, every element has a
wandering neighborhood in $\go$.
\end{proof}

\begin{corollary}
Let $G$ be a principal groupoid so that $\cs(G)$ is GCR.  The
largest Fell ideal of $\cs(G)$  is $\cs(G|_Y)$ where
\begin{equation}
Y=\{x \in \go : \text{ there exists a wandering neighborhood
of } x\}.\notag
\end{equation}
\end{corollary}

\begin{proof}
Since $G$ is principal and $\cs(G)$ is GCR, by
\lemref{idealbijection} we know every closed ideal is of the form
$\csg{Y}$
for some open G-invariant subset $Y \in \go$.  From \corref{cor1}
we see that the $Y$ defined above is G-invariant.
Also notice that $Y$ is open.  Now apply \thmref{fellthm} and we
see
that $\cs(G|_Y)$ is a Fell algebra and that any ideal
that is also a Fell algebra, must be contained in $\csg{Y}$.
\end{proof}

\thanks{This research was done as part of the
author's Ph.D. thesis at Dartmouth College under the direction of Dana~P.
Williams.  Thank you to Dana for his continued support. }

\end{document}